\documentclass[10pt,twoside]{amsart}
 \usepackage{amsmath,amssymb,amsthm,amscd,latexsym,graphicx}
 \usepackage{epigraph}
 \usepackage[OT2,OT1]{fontenc}

\usepackage{tikz}

\usepackage{graphicx}
\usepackage[all]{xy}

 \font\caps=cmcsc10                    
 \font\Caps=cmcsc10 scaled \magstep1   

 \makeatletter
 \@specialpagefalse\headheight=8.5pt
 \makeatother

 \def\TSkip{\medskip}
 \newbox\TheTitle{\obeylines\gdef\GetTitle #1
 \ShortTitle  #2
 \SubTitle    #3
 \Author      #4
 \ShortAuthor #5
 \EndTitle
 {\setbox\TheTitle=\vbox{\baselineskip=20pt\let\par=\cr\obeylines%
 \halign{\centerline{\Caps##}\cr\noalign{\medskip}\cr#1\cr}}%
         \copy\TheTitle\TSkip\TSkip%
 \def\next{#2}\ifx\next\empty\gdef\STitle{#1}\else\gdef\STitle{#2}\fi%
 \def\next{#3}\ifx\next\empty%

 \else\setbox\TheTitle=\vbox{\baselineskip=20pt\let\par=\cr\obeylines%
     \halign{\centerline{\caps##} #3\cr}}\copy\TheTitle\TSkip\TSkip\fi%
 \centerline{\caps #4}\TSkip\TSkip%
 \def\next{#5}\ifx\next\empty\gdef\SAuthor{#4}\else\gdef\SAuthor{#5}\fi%
 \catcode'015=5}}\def\Title{\obeylines\GetTitle}

 \def\Abstract{\begingroup\narrower
     \parskip=\medskipamount\parindent=0pt{\caps Abstract. }}
 
\def\Resume{\begingroup\narrower
     \parskip=\medskipamount\parindent=0pt{\caps Résumé. }}
 
 \long\def\MSC#1\EndMSC{\def\arg{#1}\ifx\arg\empty\relax\else
      {\par\narrower\noindent%
      2000 Mathematics Subject Classification: #1\par}\fi}

 \long\def\KEY#1\EndKEY{\def\arg{#1}\ifx\arg\empty\relax\else
         {\par\narrower\noindent Keywords and Phrases: #1\par}\fi\TSkip}

 \long\def\DATE#1\EndDATE{\def\arg{#1}\ifx\arg\empty\relax\else
         {\par\narrower\noindent \center{\textit{#1}}\par}\fi\TSkip\TSkip\TSkip}

 \font\bf= cmbx10 at 10pt

 \newcommand{\lra}{\longrightarrow}

 \newcommand{\F}{\mathbb{F}}
 
 \newcommand{\Bl}{\mathrm{Bl}}

 \renewcommand{\P}{\mathbb P}

 \newcommand{\Z}{\mathbb{Z}}

 \newcommand{\Stab}{{\mathrm{Stab}}}

  \newcommand{\SStab}{\mathbf{Stab}}

 \newcommand{\AAut}{\mathbf{Aut}}
 \newcommand{\Aut}{\mathrm{Aut}}

  \newcommand{\End}{\mathrm{End}}

 \newcommand{\Hom}{\mathrm{Hom}}

  \newcommand{\Ver}{\mathrm{Ver}}
   \newcommand{\Frob}{\mathrm{Frob}}

 \newcommand{\Lie}{\mathrm{Lie}}

\newcommand{\Id}{\mathrm{Id}}

 \theoremstyle{plain}
 \newtheorem{thm}{Theorem}[section]
 \newtheorem{defi}[thm]{Definition}

 \newtheorem{prop}[thm]{Proposition}

 \newtheorem{lem}[thm]{Lemma}

 \theoremstyle{remark}
 \newtheorem{rem}[thm]{Remark}

 \newenvironment{dem}{{\bf Proof.}}{\hfill$\square$}

\date{March 2010}

\newif\ifquoteopen
\catcode`\"=\active 
\DeclareRobustCommand*{"}{%
   \ifquoteopen
     \quoteopenfalse ''%
   \else
     \quoteopentrue ``%
   \fi
}
\begin{document}


 \Title
    Realisation of Abelian varieties as automorphism groups
 \ShortTitle
Abelian
\SubTitle
 \Author
  Mathieu Florence
 \ShortAuthor
 \EndTitle

\address{Mathieu Florence, Equipe de Topologie et G\'eom\'etrie Alg\'ebriques, Institut de Math\'ematiques de Jussieu,  Sorbonne Université,  Paris. }
\email{ https://webusers.imj-prg.fr/\textasciitilde mathieu.florence/}


\Abstract

Let $A$ be an Abelian variety over a field $F$. We show that $A$ is isomorphic to the automorphism group scheme of a smooth projective $F$-variety if, and only if, $\Aut_{gp}(\overline A)$ is finite. This result was proved by Lombardo and Maffei \cite{LM} in the case $F=\mathbb C$ , and recently by Blanc and Brion \cite{BB} in the case of an algebraically closed $F$.\\
\bigskip

\Resume

Soit $A$ une variété abélienne sur un corps $F$. On montre que $A$ est isomorphe au schéma en groupes des automorphismes d'une $F$-variété projective et lisse, si et seulement si le groupe des $\overline F$-automorphismes de $A$ est fini.\\
 Ce résultat est dû à Lombardo et Maffei \cite{LM} lorsque $F=\mathbb C$. Il est dû à Blanc et Brion  \cite{BB} lorsque $F=\overline F$ . \\
\tableofcontents

\section{Introduction.}
Let $F$ be a field, with algebraic closure $\overline F$. Let  $X$ be a projective variety over $F$. The automorphism group functor $\AAut(X)$ is represented by a group scheme, locally of finite type over $F$ (\cite{MO}, Theorem 3.7). Conversely, given a group scheme $G$, of finite type over $F$, it is natural to ask whether $G$ can be realised as the automorphism group of such an $X$. When $G=A$ is an abelian variety, this question was recently considered in \cite{LM}. When $F=\mathbb C$, Lombardo and Maffei prove that $A$ is the automorphism group of a projective smooth complex variety, if and only if $\Aut_{gp}(A)$ is finite. They use analytic methods. Their result was extended to $F$ algebraically closed of any characteristic in \cite{BB}, using algebro-geometric techniques: blowups, Lie algebra computations and modding out actions of finite group schemes. Making a different use of these tools, we provide  a generalisation of this result, to the case of all ground fields $F$.\\
\subsection{Sketch of our construction. }\hfill\\
Let $A/F$ be an abelian variety over a field $F$, such that $G:=\Aut_{gp}(\overline A)$ is finite. We first introduce an integer $n \geq 1$, invertible in $F$, such that $G$ acts faithfully on the $n$-torsion subgroup $A[n](\overline F)$. \\Then, we pick an abelian variety $B_1/F$, enjoying the following properties.
\begin{enumerate}
    \item {The abelian varieties $A$ and $B_1$ are `orthogonal', in the sense that $$\Hom_{gp}(\overline A, \overline B_1)=0,$$ where homomorphisms are taken over $\overline F$.} \item { There exists an injection (of algebraic $F$-groups) $$\iota: A[n] \hookrightarrow B_1. $$
}
\end{enumerate}
We denote by $B_2/F$ the abelian variety fitting into the diagonal extension \[ 0 \lra A[n] \stackrel {a \mapsto (a,\iota(a))} \lra A \times B_1 \stackrel \pi \lra B_2 \lra 0. \]
Using point (1) above, we prove that automorphisms of (the variety) $B_2$ are diagonal: they come from automorphisms of $ A \times B_1 $, respecting orbits under the embedded $A[n]$.
Next, we build an appropriate smooth closed $F$-subvariety $Y_2 \subset B_2$, stable by translations by $A \simeq \pi(A\times \{0\}) \subset B_2$.\\ We define a smooth $F$-variety $X$ as the blowup \[X:=\Bl_{Y_2} B_2.\]
The natural arrow \[  A \lra \AAut(B_2),\]  given by translations,  lifts to an arrow \[ \tau: A \lra \AAut(X).\] We show that $\tau$ is an isomorphism of algebraic groups over $F$.

\section{Notation.}
\subsection{Geometry over $F$.}
Let $F$ be a field, with algebraic closure $\overline F$, and separable closure $F_s \subset\overline F$. We denote by $F[\epsilon]$, $\epsilon^2=0$, the $F$-algebra of dual numbers. We use it for differential calculus.\\By a variety over $F$, we mean a separated $F$-scheme of finite type. \\An algebraic $F$-group (or simply $F$-group) is an $F$-group scheme of finite type. It is often assumed to be reduced, hence smooth over $F$.\\
Let $X$ be a variety over $F$. For a field extension $E/F$, we denote by $X_E:=X \times_F E$ the $E$-variety obtained from $X$ by extending scalars. We put $\overline X:=X\times_F \overline F$. \\If $X$ is smooth over $F$, we denote by $TX \lra X$ the tangent bundle of $X$. A global section of the tangent bundle  is called a vector field on $X$. \\ We denote by $\Aut(X)$ the (abstract) group of automorphisms of the $F$-variety $X$, and by ${\Aut}(\overline X)$ the group of automorphisms of the $\overline F$-variety $\overline X$.  \\
If $X/F$ is a projective variety, we denote by $\AAut(X)$ the  $F$-group scheme of automorphisms of $X$; it is locally of finite type over $F$. By \cite{B}, Lemma 3.1, there is a canonical isomorphism \[ H^0(X,TX) \stackrel \sim \lra \Lie(\AAut(X)).\] \\
If an abstract group $G$ acts on a variety $X$, and if $Z \subset X$ is a closed subvariety, we denote by $\Stab_G(Z) \subset G$, or simply by $\Stab(Z) \subset G$ when no confusion arises, the subgroup of transformations leaving $Z$ (globally) invariant.\\ Let $G/F$ be a group scheme, locally of finite type. In the situation where $G$ acts on $X$, we use the notation  $\SStab_G(Z) \subset G$ for the  closed $F$-subgroup scheme defined by \[\SStab_G(Z) (A)=\{g \in G(A), g(Z_A)=Z_A\}, \]  for all commutative $F$-algebras $A$. That it is representable follows from \cite{DG}, II 1.3.6.
\subsection{Frobenius and Verschiebung.}
If $F$ has characteristic $p>0$, we put $$X^{(1)}:=X \times_{\Frob} F,$$  extension of scalars taken with respect to $\Frob: F\stackrel {x \mapsto x^{p}} \lra F.$ \\ Recall the Frobenius homomorphism $$\Frob_X : X \lra X^{(1)};$$ it is a morphism of $F$-varieties, functorial in  $X$. \\If  $X/F$ is an algebraic group, it is a group homomorphism.\\ If $X$ is a commutative algebraic group, there is the Verschiebung homomorphism  $$\Ver_X : X^{(1)} \lra X,$$ satisfying $(\Ver_X \circ \Frob_X)=p\Id_X$. \\If moreover $X/F$ is a semi-abelian variety, $\Ver_X$ and $\Frob_X$ are isogenies.
\subsection{Abelian varieties.}
If $A$  and $B$ are Abelian varieties over $F$, we denote by $\Hom_{gp}(A,B)$ the group of homomorphisms of algebraic $F$-groups, from $A$ to $B$. We denote by ${\Hom_{gp}}( \overline A,\overline B)$  the group of homomorphisms of algebraic $\overline F$-groups, from $\overline A$ to $\overline B$. These are finite free $\Z$-modules. We adopt the similar notation for endomorphisms ($\End_{gp}$) and automorphisms ($\Aut_{gp}$). For an integer $n \geq 1$, we denote by $A[n]$ the $n$-torsion of $A$, seen as a finite group scheme over $F$.\\
\subsection{Barycentric operations.}
Let $A$ be an abelian variety over $F$. Then $A$ comes naturally equipped with barycentric operations with integer coefficients. More precisely, for a positive integer $n$, denote by $$\mathbf Z^n_1 \subset \mathbf Z^n$$ the subset consisting of integers $\alpha=(\alpha_1,\ldots, \alpha_n)$, with $\alpha_1+\ldots+\alpha_n=1$.\\ For $\alpha \in \mathbf Z^n_1$, there is a barycentric operation \[ \mathcal B_\alpha: A^n \lra A,\] \[ (x_1,\ldots,x_n) \mapsto \alpha_1 x_1+\ldots+\alpha_n x_n.\]  Associativity  of the group law of $A$, provides natural associativity relations  between the $\mathcal B_\alpha$'s, for various $\alpha's$.\\ For instance,  pick $\alpha=(\alpha_1, \alpha_2) \in \mathbf Z^2_1$ and  $\gamma=(\gamma_1,\gamma_2) \in \mathbf Z^2_1$,   and set $$\delta:=(\alpha_1 \gamma_1,\alpha_2 \gamma_1,\gamma_2) \in \mathbf Z^3_1.$$
Then, we have the associativity rule
\[\mathcal B_\gamma (\mathcal B_\alpha (x_1,x_2),x_3)= \mathcal B_\delta (x_1,x_2,x_3).\] 

\begin{rem}
More generally, these barycentric operations exist for \textit{torsors} under commutative algebraic $F$-groups.
\end{rem}

\begin{defi}
Let  $X \subset A$  be an $F$-subvariety. We say that $X$ is stable under all barycentric operations, if the restriction $$(\mathcal B_\alpha)_{\vert X^n}:X^n \to A$$ factors through the closed immersion $X \hookrightarrow A$, for every $n\geq 2$ and every  $\alpha \in \mathbf Z^n_1$.\\ In this case, we also say that $X$ is \textit{barycentric}. 
\end{defi}
Note that $X$ is barycentric if and only if it is a translate of an algebraic $F$-subgroup $ \overrightarrow X \subset A$.   Checking this fact is left as an exercise for the reader. Of course, $X(F)$ might be empty. If $X$ is geometrically reduced and geometrically connected, so is $ \overrightarrow X$- hence $ \overrightarrow X$ is an abelian subvariety of $A$.\\
Let $A$ and $B$ be two abelian varieties over $F$. Recall the essential fact $$\Hom_{\overline F-var}(\overline A,\overline B)= B(\overline F) \times \Hom_{gp}(\overline A,\overline B).$$In particular, morphisms (of varieties) between abelian varieties commute with the barycentric operations $\mathcal B_\alpha$. \\ If $X \subset A$ is a geometrically reduced closed $F$-subvariety, the smallest geometrically reduced barycentric $F$-subvariety containing $X$ is called the barycentric envelope of $X$. We denote it by $\mathcal E(X)$. \\
Assume now that $X$ is geometrically reduced and geometrically connected. Pick $n \geq 1$ and $\alpha \in \mathbf Z^n_1$. Consider $\mathcal B_\alpha(X^n)  \subset A$ as a geometrically reduced and geometrically connected closed subvariety of $A$. Then, if $n$ and $\alpha$ are chosen so that $\mathcal B_\alpha(X^n)$ is of maximal  dimension, we have $\mathcal B_\alpha(X^n)= \mathcal E(X)$. Thus, $ \mathcal E(X)$, being geometrically connected and geometrically reduced, is a translate of an abelian subvariety of $A$.

\section{Statement of the theorem.}
\begin{thm}
Let $A$ be an Abelian variety, over a field $F$.  The following are equivalent:\\
1) The group $G:={\Aut_{gp}}(\overline A)$ is finite.\\
2) There exists a smooth projective $F$-variety $X$, such that $A$ is isomorphic to $\AAut(X)$ (as algebraic groups over $F$).
\end{thm}

Note that $2) \Rightarrow 1)$ can be checked over $\overline F$, which follows from \cite{BB}, Theorem A.\\ Our task in this paper is to prove the converse implication.
\section{Auxiliary results.}
\subsection{Blowups.}

This section contains two elementary lemmas on automorphisms of blowups, which we provide with short proofs. A good recent reference on this topic, also containing more advanced material, is section 2 of  \cite{M}.
\begin{lem}\label{blowup}\label{AutBlow}
Let $Y \hookrightarrow D$ be a closed immersion of smooth $F$-varieties, such that all connected components of $Y$ have codimension $\geq 2$ in $D$. \\ Denote by $\beta: X:=\Bl_Y(D) \lra D$ the blowup of $Y$ inside $D$. \\The $F$-variety $X$ is smooth.\\Let $f$ be an automorphism of the $F$-variety $D$. Then, $f$ lifts via  $\beta$ to an automorphism of $X$, if and only if  $f(Y)=Y$.\\

\end{lem}
\begin{dem}
If $f(Y)=Y$, then $f$ lifts to an automorphism of $X$ by the universal property of the blowup.\\
Conversely, assume that $f$ lifts to an automorphism $\phi$ of $X$, so that we have a commutative square\[ \xymatrix{  X \ar[r]^\phi \ar[d] &  X \ar[d] \\   D \ar[r]^f &  D.}\] To check  that  $f(Y)=Y$, can assume that $F=\overline F$. It then suffices to prove that $Y \subset D$ and $f(Y)  \subset D$ have the same set of $F$-rational points. This is clear, since the fiber of $\beta$ over a point $x \in D(F)$  is either a point if $s \notin Y(F)$, or a projective space of dimension $\geq 1$ if $x \in Y(F)$.
\end{dem}

Lemma \ref{AutBlow} has an infinitesimal analogue, as follows.
\begin{lem}\label{blowup}
Let $Y \hookrightarrow D$ be a closed immersion of smooth $F$-varieties, such that all connected components of $Y$ have codimension $\geq 2$ in $D$. \\Denote by $\beta:  X:=\Bl_Y(D) \lra D$ the blowup of $Y$ inside $D$. \\Let $s: D \lra TD$ be a vector field on $D$. Then, $s$ lifts to a vector field on $X$, if and only if $s_{\vert Y} $ takes values in $TY$. 
\end{lem}
\begin{dem}
Denote by $i: E \hookrightarrow X$ the exceptional divisor. The restriction $$\beta_{\vert X-E}: X-E \lra D-Y$$ is an isomorphism.\\
We thus have a natural injective $F$-linear arrow \[ \rho: H^0(X,TX) \lra H^0(D-Y,TD) =H^0(D,TD),\] \[ \sigma \mapsto \sigma_{\vert X-E}.\] Note that the equality $H^0(D-Y,TD) =H^0(D,TD)$ follows from the fact that $Y \subset D$ has codimension $\geq 2$. On $E$, we have a natural extension of vector bundles \[ 0 \lra TE \lra i^*(TX) \lra N_{E/X} \lra 0,\] where $N_{E/X}  \simeq \mathcal O_E(-1)$ is the normal bundle of $E$ in $X$.  Since $Y$ has codimension $\geq 2$ in $D$, we have $H^0(E, O_E(-1))=0$. This can be checked on the fibers of $\beta$ over geometric points of $Y$, which are projective spaces of dimension $\geq 1$. Hence, $\sigma_{\vert E}$ takes values in $TE$. Consequently, $\rho(\sigma)_{\vert Y}$ takes values in $TY$.\\
Conversely, let $s: D \lra TD$ be a vector field on $D$. Then $s$ corresponds to an automorphism $\psi$ of the $F[\epsilon]$-scheme $D \times_F F[\epsilon]$, reducing to the identity at $\epsilon=0$. Assume that $s_{\vert Y} $ takes values in $TY$. Then, $\psi$ restricts to an automorphism of the closed subscheme $Z \times_F F[\epsilon] \subset D \times_F F[\epsilon]$. By the universal property (and compatibility with base change) of the blowup, $\psi$ lifts, via $\beta \times_F F[\epsilon]$, to an automorphism of  $X \times_F F[\epsilon]$. Equivalenty, $s$ lifts, via $\beta$, to a vector field on $X$.

\end{dem}
\subsection{Hypersections in projective space.}\hfill \\
We could not find a reference in the literature for the following result, so that we provide it with a proof.
\begin{prop}\label{exo}
Let $S$ be a geometrically irreductible smooth projective $F$-variety, of dimension $\geq 2$. Let $m \geq 1$ be an integer. Then, $S$ contains a   geometrically irreductible smooth projective $F$-curve, of genus $g\geq m$.
\end{prop}
\begin{dem}
 Pick a projective embedding $S \subset \P^n$ (everything is over $F$). Let $d \geq 1$ be an integer. Let $H \subset  \P^n$ be a degree $d$ hypersurface, given by $h \in H^0(\P^n,\mathcal O(d))$.  By Bertini's theorem, for $d$ large enough and $h$ general, $S \cap H$ is smooth and  geometrically irreductible, of dimension one less than $S$. This version of Bertini's theorem works over any $F$- see \cite{P} and \cite{CP} for the delicate case where $F$ is finite. Proceeding by induction, we reduce to the case where $S$ is a surface.\\ We then take $C:=S \cap H$,  and show that $g(=h^1(C,\mathcal O_C))$ goes to infinity with $d$.
To do so, consider the exact sequence of coherent $\mathcal O_{\P^n}$-modules \[0 \to \mathcal O_S(-d) \xrightarrow{\times h}  \mathcal O_S \to \mathcal O_C \to 0.\] Taking Euler characteristics, we get \[g-1=-\chi(  \mathcal O_C)= \chi( \mathcal O_S(-d))-\chi(\mathcal O_S).\] We conclude using the following fact, applied to $X=S$.\\ For a closed $m$-dimensional $F$-subvariety $X  \subset \P^n$, the association \[d \mapsto \chi( \mathcal O_X(-d)) \] is a degree $m$ polynomial function of $d$. A classical proof is by induction on $m\geq 0$.
\end{dem}

\subsection{(Semi-)abelian varieties.}
The next Lemma is borrowed from \cite{B}, Lemma 5.3. We provide here  a different proof. In practice, we will apply it to abelian varieties, in which case it is due to Chow.
\begin{lem}\label{separable}
Assume that $F$ has characteristic $p>0$.\\
Let $A,B$  be  semi-abelian varieties over $F$. Then, all elements of ${\Hom_{gp}(\overline A,\overline B)}$ are defined over the separable closure $F_s \subset \overline F$.
\end{lem}

\begin{dem}
We have to show the following. Let $E/F$ be a purely inseparable algebraic extension. Let $g: A_E \lra  B_E$ be a homomorphism of algebraic groups over $E$. Then $g$ is defined over $F$. Without loss of generality, we can assume that $E/F$ is finite. By induction, we reduce to the case where $E=F(\sqrt[p] a) /F$ is a primitive purely inseparable extension of height one.  Note that $\Frob: E \lra E$ takes values in $F$. Hence,  $g^{(1)}: A_E^{(1)} \lra B_E^{(1)}$ is defined over $F$. The Frobenius homomorphism $$\Frob_A: A \lra  A^{(1)}$$  presents $A^{(1)}$  as a quotient of $A$, by a finite (characteristic) sub-$F$-group $\mu_{A} \subset A$.  \\From the relation $$\Ver_A \circ \Frob_A =p \Id_A,$$ we deduce  $\mu_{A} \subset A[p]$. Same holds  for $B$.\\ Combining these facts, we get that the  $E$-morphism $$ A/\mu_{A} \lra B/\mu_{B}$$ induced by $g$, is defined over $F$. Modding out further, we get that the $E$-morphism $$ A/A[p] \lra  B/B[p],$$ induced by $g$, is defined over $F$. \\Via the iso $$ A/A[p] \stackrel \sim \lra A$$ $$
\overline{a} \mapsto pa,$$ this isomorphism is actually $g$ itself. The Lemma is proved.
\end{dem}
\begin{lem}\label{Mustbetrivial}
For each $n\geq 2$, there exists an (absolutely) simple $n$-dimensional abelian variety $A$ over $F_s$.
\end{lem}
\begin{dem}
Since $F_s$ is separably closed, `simple' is the same as `absolutely simple', for abelian varieties over $F_s$ (use Lemma \ref{separable}). Without loss of generality, we assume that $F_s$ is the algebraic closure of its prime subfield.  Over $\overline {\mathbb Q}$, we can then  use the existence of abelian surfaces with a prescribed CM type. Over $\overline {\F}_p$, we can use Honda-Tate theory. For concrete constructions, and more general results,  we refer to \cite{M}, Theorem 1 (where $F_s=\overline {\mathbb Q}$), and \cite{HZ},  Theorem 2 (where $F_s=\overline {\F}_p$).

\end{dem}
\begin{lem}\label{smooth}
Let $B$  be an abelian variety over $F$, whose simple factors (over $\overline F$) are of dimensions $\geq 2$. (Equivalently: all $\overline F$-homomorphisms from an elliptic curve to $\overline B$ are constant.) \\Then, there exists a smooth $F$-subvariety $Y \subset B$, which is a disjoint union of smooth $F$-curves, and of a separable closed point, such that $$\SStab(Y)=\{\Id\} \subset{\AAut(  B)}.$$
\end{lem}
\begin{dem}
Assume first that $B$ is $F$-simple, in the sense that it has no non-trivial proper abelian $F$-subvariety.
By Proposition \ref{exo}, we can pick a geometrically irreducible smooth $F$-curve  $C \subset B$, of arbitrarily large genus $g\geq 2$.\\ The group $ {\Aut(\overline C)}$ is finite.  Indeed, $\Lie(\AAut(C))$ is the space of vector fields on $C$, which vanishes since $g \geq 2$.\\ Let us show that  $\mathcal E(C)=B$. The barycentric envelope $\mathcal E(C)$ is a translate of an abelian subvariety $B' \subset B$. Since $B$ is $F$-simple, we get $B'= B$, hence $\mathcal E(C)=B$. \\ Now, let $g \in  \AAut ( B)(\overline F[\epsilon])=B(\overline F[\epsilon]) \times \Aut_{gp}(\overline B)$ be such that $$g_{\vert C \times_F \overline F[\epsilon] }=\Id_{\vert C \times_F \overline F[\epsilon]}.$$ Because $g$ commutes to barycentric operations, $g$ acts as the identity on the closed subscheme $$\mathcal E(C) \times_F  \overline F[\epsilon] \subset B \times_F  \overline F[\epsilon].$$ Since $\mathcal E(C)=B$, it follows that   $g=\Id$. Thus, we get a natural embedding of $F$-group schemes $$H:=\SStab_{{ \AAut (B)}}(C)  \hookrightarrow \AAut(C).$$ In particular, $H$ is finite étale over $F$. Let $E/F$ be a finite separable field extension, such that $H(E)=H(\overline F)$.  Denote by $$\Phi:=\bigcup_{h \in H(E), h \neq e} \overline B ^h \subset \overline B$$ be the (strict) closed subscheme, consisting of points fixed by at least one non-trivial element $h \in H(E)$. It is defined over $F$ by Galois descent. There exists a finite separable field extension $L/E$, and a point $b \neq 0 \in B(L)$, which does not lie in $\Phi(L)$, nor in $C(L)$. We then have a separable zero-cycle $[b]$ in the $F$-variety $B$, of degree $[L:F]$. Define $Y \subset B$ as the disjoint union of  $[b]$ and $C$. We claim that $Y$ has the required property. Indeed, let $f \in  {\AAut}(B)(\overline F[\epsilon])$ be an automorphism stabilizing $Y$- or more accurately, $Y \times_F \overline F[\epsilon] \subset B \times_F  \overline F[\epsilon] $. Then, $f$ permutes the two connected components of the scheme $Y \times_F \overline F[\epsilon]$. For dimension reasons, it preserves $C \times_F \overline F[\epsilon]$ on the one hand, and  $[b] \times_F \overline F[\epsilon]$ on the other hand. From the first fact, we know that $f$  belongs to $H(\overline F)$; in particular, it is defined over $E$, hence over $L$. From the latter fact, we get $f(b)=b$, hence $f=\Id$. The Lemma is proved in this case.\\
 Assume now that $B=B_1 \times \ldots B_n$, where the $B_i$'s are $F$-simple abelian varieties. We can then adapt the preceding proof, as follows. For each $i$, let $C_i \subset B_i$, $L_i/E_i/F$ and $b_i \in B (L_i)$ be as in the first part of the proof. We can fulfill the extra requirements that no $C_i$ passes through $0$, and that the  $C_i$'s are of different genus (using Proposition \ref{exo}). In particular, when $i \neq j$,  $\overline C_i$ is not  $\overline F$-isomorphic to  $\overline C_j$. We can also assume that $L_i=L$ and $E_i=E$ are independent of $i$. Set \[b:=(b_1,\ldots,b_n) \in B(L).\] Define $Y$ to be the disjoint union of $[b]$, and of the $n$ curves $$C_i \simeq \{0\} \times \ldots \times \{0\} \times C_i  \times \{0\} \times \ldots \times \{0\}  \hookrightarrow B.$$ It is not hard to see, that $Y$ enjoys the required property.\\
In general, write $B=(\prod_1^r B_j)/\mu$, where $S_1,\ldots,S_r$ are $F$-simple abelian varieties, and where $\mu$ is a finite $F$-subgroup, intersecting trivially each coordinate axis. We  can choose $$ Y \hookrightarrow B_1 \times \ldots B_n$$ as in the previous part of the proof, and such that the composite $$Y \hookrightarrow B_1 \times \ldots B_n \stackrel{can} \lra  (\prod_1^r B_j)/\mu =B$$ is a closed immersion, identifying $Y$ to a smooth closed subvariety of $B$.\\ An automorphism of $B$ stabilizing  $Y \subset B$ then lifts, via the quotient $can$, to an automorphism of $B_1 \times \ldots B_n$  stabilizing $Y \subset B_1 \times \ldots B_n$. We conclude as before.
\end{dem}

\section{Proof of the implication $1) \Rightarrow 2)$.}
Let $A/F$ be an abelian variety, such that $G:=\Aut(\overline A)$ is finite. We give a construction of a smooth projective $F$-variety $X$, such that $A=\AAut(X)$, in several steps.
\subsection{Construction of $X$.}
\hfill\\
Denote by $g$ the dimension of $A$.\\
Let $n\geq1$ be an integer, invertible in $F$, such that the action of $G$ on  $A[n](F_s)\simeq (\Z/n)^{2g}$ is faithful.  Such an $n$ exists: use that $G$ is finite, and that torsion points of  order prime to $\mathrm{char}(F)$ in $A(\overline F)$ are Zariski-dense in $A$.\\
Let $B_s$ be an abelian variety over $F_s$,  of dimension $g' \geq g$, such that $$ {Hom}_{gp}(\overline A,\overline B)= {Hom}_{gp}(\overline B,\overline A)=0.$$ Since $\overline A$  has a finite number of simple components (up to isogeny), which are all defined over $F_s$ by Lemma \ref{separable}, the existence of $B_s$ follows from Lemma \ref{Mustbetrivial}. For example, take for $B_s$ a product of simple abelian varieties, of dimensions greater than that of the simple components of $\overline A$.\\

Let $E/F$, be the finite Galois extension, with group $\Gamma$, which is minimal w.r.t. the following properties.

\begin{enumerate}
    \item{The extension $E/F$ splits the $F$-group of multiplicative type $A[n]$.\\ In other words, $A[n](E)\simeq (\Z/n)^{2g}$.} \item{ The abelian variety $B_s$ is defined over $E$: there exists an abelian $E$-variety $B_E$, such that $B_E \times_E F_s \simeq B_s$.} \item{Same as (1), for $B_E$: we have  $B_E[n](E)\simeq (\Z/n)^{2g'}$.}
\end{enumerate}

Using (1), we view  $A[n](E)$ as a $(\Z/n) [\Gamma]$-module.\\
Introduce  the Weil restriction of scalars $$B_1:=R_{E/F}(B_E).$$   Geometrically, we have $\overline B_1 \simeq \overline B_s^{m}$, where $m$ is the cardinality of $\Gamma$. \\ We have $$B_1 [n]= R_{E/F}((\Z/n)^{2g'}),$$ so that $E/F$ splits  $B_1 [n]$, and $B_1 [n](E)$ is a free $(\Z/n ) [\Gamma]$-module of rank $2g'$.\\
\begin{lem}
There exists an embedding of $(\Z/n \Z) [\Gamma]$-modules $$A[n ](E) \hookrightarrow B_1 [n](E);$$ that is to say,  an embedding of finite \'etale $F$-group  schemes $$\iota: A[n]\hookrightarrow   B_1 [n].$$
\end{lem}
\begin{dem} We give two (seemingly) different proofs.\\
The first one uses the perfect duality $$(.)^\vee:=\Hom(., \Z/n),$$ in the category of  $(\Z/n ) [\Gamma]$-modules. Pick a generating set  $t_1,\ldots, t_{2g'}$ of the $\Z/n$-module $A[n ](E) ^\vee$- which is free of rank  $2g \leq 2g'$. Introduce the surjection of $(\Z/n ) [\Gamma]$-modules \[(\Z/n ) [\Gamma]^{2g'} \lra  A[n ](E) ^\vee,\] \[ e_i \mapsto t_i,\] where $e_i$ denotes the $i$-th element of the canonical basis. Dualizing it yields an injection of $(\Z/n ) [\Gamma]$-modules \[\iota: A[n ](E) \lra  (\Z/n ) [\Gamma]^{2g'} \simeq B_1 [n],\] concluding the construction. \\
The second proof is more conceptual. Choose an embedding of constant $E$-group schemes \[(\Z /n)^{2g} \simeq  A_E[n] \hookrightarrow  B_E[n] \simeq (\Z /n)^{2g'},\] which exists simply because $g \leq g'$. \\Applying $R_{E/F}$ yields an embedding of $F$-group schemes \[  R_{E/F}(A_E)[n] \hookrightarrow  R_{E/F}(B_E)[n]=B_1[n].\] Composing it with the natural embedding of $F$-groups \[ A[n] \hookrightarrow R_{E/F}(A_E)[n], \] arising by adjunction from the identity of $A_E[n]$, we get the desired $\iota$.
\end{dem}

Form the exact sequence of algebraic $F$-groups \[ 0 \lra A[n] \stackrel {a \mapsto (a,\iota(a))}\lra A \times B_1  \stackrel \pi \lra B_2 \lra 0.\] Its cokernel $B_2 $ is an abelian variety over $F$. \\We have $F$-embeddings $$A \stackrel {a \mapsto (a,0)} \hookrightarrow B_2$$ and $$B_1 \stackrel {b_1 \mapsto (0,b_1)} \hookrightarrow B_2.$$
Introduce the quotient $$q: B_2 \lra B_3:=B_2/A \simeq B_1/\iota(A[n]).$$
Let $Y_3 \subset B_3$ be a smooth $F$-subvariety, enjoying the properties of Lemma \ref{smooth}, where we take $B$ to be our $B_3$, and set $Y_3:=Y$.\\  Put $$Y_2:=q^{-1}(Y_3).$$ The restriction $$q_{\vert Y_2 }: Y_2 \lra Y_3$$ is an $A$-torsor.\\ We now define $$X:=\Bl_{Y_2}(B_2)$$ to be the blowup of $Y_2$ in $B_2$. 

\subsection{Proof that $\AAut(X) \simeq A$.}
\hfill\\
Translating by elements of $A$ inside $B_2$ yields a natural arrow $$ A \lra \AAut(B_2). $$ Since $Y_2\subset B_2$ is stable by these translations,  we get an induced  arrow of  $F$-group schemes $$\tau: A \lra \AAut(X). $$ It is clear that $\tau$ is an embedding. We are going to show that it is an isomorphism.\\ Let us first check that it induces a bijection $$A(\overline F) \stackrel \sim \lra  {\Aut}( \overline X)=\AAut(X)(\overline F).$$ Pick $\phi \in  {\Aut}( \overline X)$. It induces a birational isomorphism $f_2$ of the $\overline F$-variety $\overline B_2$, which is a regular isomorphism since $B_2$ is an abelian variety. Thus, we get a commutative diagram \[ \xymatrix{ \overline X \ar[r]^\phi \ar[d] & \overline  X \ar[d] \\ \overline  B_2 \ar[r]^{f_2} & \overline  B_2,}\] where the vertical arrows are the structure morphism of the blowup.\\  Using Lemma \ref{AutBlow}, we get $f_2( \overline Y_2)= \overline Y_2$.  We know that $$f_2(x)=g_2(x)+t_2,$$ where $g \in  {\Aut}_{gp}( \overline B_2)$, and $t_2 \in B_2(\overline F)$. We have to show that $g_2=\Id$  and $t_2 \in A(\overline F)$. To do so, we can assume without loss of generality that $t_2 \in B_1(\overline F)$.\\ We then have to prove $g_2=\Id$ and $$t_2\in A(\overline F) \cap B_1(\overline F)=\iota(A[n])(\overline F).$$\\
Geometrically, $ \overline B_1 \simeq  \overline {B}_s^{m}$. Since  $ {Hom}_{gp}(\overline A,\overline {B}_s)= {Hom}_{gp}(\overline {B}_s,\overline A)=0$, we get $$ {Hom}_{gp}(\overline A,\overline B_1)= {Hom}_{gp}(\overline B_1,\overline A)=0.$$ Therefore $g_2$ leaves $\overline A \subset \overline B_2$ and $\overline B_1 \subset \overline B_2$ stable. \\We infer that $g_2$ lifts, via $\overline \pi$, to  a diagonal group automorphism  $$\delta=(h,g_1)$$ of $\overline A \times \overline B_1$,
which automatically leaves the diagonally embedded $\overline A[n]$ stable. \\Consider the automorphism of $\overline B_1$ given by  $$f_1(b_1):= g_1(b_1)+t_2,$$ and the diagonal automorphism of $\overline A \times \overline B_1$ given by $$ \Delta(a,b_1):=(h(a), f_1(b_1)).$$ Since $\delta$ leaves  $\overline A \times \iota(\overline A[n]) \subset \overline A \times \overline B_1$ stable, there exists $f_3\in \Aut(\overline B_3)$ such that the diagram \[ \xymatrix{ \overline A \times  \overline B_1 \ar[r]^{\Delta} \ar[d]^{\overline \pi} &  \overline A \times  \overline B_1 \ar[d]^{\overline \pi} \\ \overline  B_2 \ar[r]^{f_2} \ar[d]^{\overline q}  & \overline  B_2 \ar[d]^{\overline q} \\ \overline B_3 \ar[r]^{f_3}  & \overline B_3 }\]  commutes.

Because $f_2(\overline Y_2)=\overline Y_2$, we get $f_3(\overline Y_3)=\overline Y_3$. By Lemma \ref{smooth}, we conclude that $f_3=\Id$. Hence, we have $t_2\in \iota(A[n])(\overline F)$ and $g_2=\Id$. Since $\delta$ preserves the diagonally embedded $\overline A[n]$, we get that $h$, restricted to $\overline A[n] \subset \overline A$, is the identity. Since $G$ acts faithfully on $A[n]$, we conclude that $h=\Id$. Hence, $g_2=\Id$ as well, and our job is done.\\

We have proved that $\tau$ induces a bijection on $\overline F$-points. If $F$ has characteristic zero, this is enough to conclude that $\tau$ is an isomorphism of algebraic $F$-groups. In general, it remains to check that the $F$-linear map on tangent spaces $$d_e(\tau): \Lie(A) \lra \Lie(\AAut(X)) $$ is bijective. Recall that $\Lie(\AAut(X)$ is the space of vector fields on $X$; that is, global section of the tangent bundle $TX \lra X$. Let $$s: X \lra TX$$ be such a section. Restricting $s$ to the complement of the exceptional divisor, we get a global section $\sigma'$ of the tangent bundle of $B_2-Y_2$. Since $B_2$ is an abelian variety, its tangent bundle is trivial, so that $\sigma'$ is given by an arrow of $F$-varieties $$\sigma': B_2-Y_2 \lra \mathbb A(\Lie(B_2)),$$ with target an affine space of dimension $\dim(B_2)$. Since $Z_2$ has codimension $\geq 2$ in $B_2$, $\sigma'$ extends to a morphism $$\sigma: B_2 \lra \mathbb A(\Lie(B_2)),$$ which is constant because $B_2/F$ is proper. Write $\sigma=t$, with $t \in \Lie(B_2)$. To conclude, we have to show $t \in \Lie(A)$. \\For $y \in B_2(\overline F)$, denote by $$ \alpha_y: \overline B_2 \lra \overline B_2$$ the $\overline F$-morphism given by $$ x\mapsto x+y.$$ Recall that the linear isomorphisms $$d_y\alpha_{-y}: T_y(\overline B_2)\stackrel \sim  \lra \Lie(B_2) \otimes_F \overline F $$ are used to trivialize the tangent bundle of $B_2$.\\  Since $\sigma$ lifts to a section of the tangent bundle of the blowup $\Bl_{Y_2}(B_2)$, Lemma \ref{blowup} implies, when $y \in Y_2(\overline F)$, that $t$ belongs to  $$d_y\alpha_{-y}(T_y(\overline Y_2)) \subset  \Lie(B_2) \otimes_F \overline F.$$  Taking a $y$ lying above (via $\overline q$) an isolated separable point of $\overline{Y}_3$,  we conclude that $t\in \Lie(A)$, as desired.

\section{Acknowledgments.}
We are grateful to the referee for  meaningful comments, leading to an improved exposition.\\
We thank Michel Brion for his careful reading, and for helpful suggestions concerning blowups. We thank Frans Oort and Olivier Wittenberg for helping us understand why Lemma \ref{Mustbetrivial} is true, and Daniel Bertrand for pointing out an existing reference for its proof.

\end{document}